\documentclass[10pt,11pt]{article}
\usepackage{amssymb}

\usepackage[dvips]{graphicx}
\usepackage{amsmath}

\providecommand{\U}[1]{\protect\rule{.1in}{.1in}}
\providecommand{\U}[1]{\protect\rule{.1in}{.1in}}
\newtheorem{theorem}{Th\'eor\`eme}[subsection]

\newtheorem{corollary}[theorem]{Corollaire}

\newtheorem{lemma}[theorem]{Lemma}

\newtheorem{remark}[theorem]{Remarque}
\newtheorem{fait}[theorem]{Fait}

\def \R{\mathbb{R}}

\def \|{\Lambda}

\def \cH{\mathcal H}

\def \fdem{$\Box$}
\def \|{{\mathcal{L}}}

\begin{document}

\begin{center}
{\Large \bf  Autour de l'exposant critique d'un groupe kleinien}

\vspace{2mm}

 Marc Peign\'e$^($\footnote{Marc Peign\'e, LMPT, UMR 6083, Facult\'e des Sciences etTechniques, Parc de Grandmont, 37200 Tours -- mail : peigne@lmpt.univ-tours.fr}$^)$,
 
\vspace{2mm}

 Journ\'ees Platon, 11, 12 \& 13 Octobre 2010

\end{center}

{\it 
Nous exposons ici quelques propri\'et\'es \'el\'ementaires autour des exposants critiques des groupes discrets d'isom\'etries en courbure strictement n\'egative et pinc\'ee. Nous rappelons quelques r\'esultats essentiels et  pr\'esentons quelques  outils classiques utilis\'es sur ce sujet}.

\section{Exposant critique et s\'erie de Poincar\'e}

Nous consid\'erons une vari\'et\'e compl\`ete simplement connexe $X$  \`a courbures sectionnelles pinc\'ees $-b^2\leq K\leq -a^2<0$ (et plus g\'en\'eralement un espace m\'etrique CAT(-1) propre). Nous fixons une origine ${o}$ de $X$ ; la distance entre deux points $x$ et $ y$ de  $X$ est not\'ee $(x, y)$. Enfin, $\partial X$ d\'esigne le bord visuel de $X$ et l'on pose $\bar X:= X\cup \partial X$.

Nous consid\'erons un groupe Kleinien de $X$  c'est-\`a-dire un groupe discret  $\Gamma$ d'isom\'etries de $X$ pr\'eservant l'orientation. Nous supposerons pour simplifier que $\Gamma$ est sans torsion : ses \'el\'ements sont donc, soit des isom\'etries paraboliques, soit des isom\'etries hyperboliques.

L'orbite $\Gamma\cdot {o}$ s'accumule sur $\partial X$ ; l'ensemble de ses valeurs d'adh\'erence est appel\'e {\it ensemble limite de $\Gamma$} et est not\'e  $\Lambda_\Gamma$.
 Cet ensemble contient $1, 2$ ou une infinit\'e non d\'enombrable de points ; dans le premier cas, $\Gamma$ est parabolique, dans le second cas il est cyclique et engendr\'e par une isom\'etrie hyperbolique et dans le dernier cas, il est dit {\it non \'el\'ementaire} et il contient en particulier des sous-groupes libres du type $\mathbb F^2$ (par exemple des sous-groupes de Schottky).
 
On \'etudie le comportement de la fonction orbitale de $\Gamma$ d\'efinie par :
 $$
 \forall {x, y}\in X, \forall R>0 \quad N_\Gamma({x}, {y}, R):= \sharp
 \{ \gamma \in \Gamma : ({x}, \gamma\cdot{y})\leq R\} 
 $$
(on \'ecrira pour simplifier $N_\Gamma(R):= N_\Gamma({o}, {o}, R)$). L'exposant $ \displaystyle \limsup_{R\to +\infty}{1\over R}\ln N_\Gamma({x}, {y}, R)$ ne d\'epend ni de ${x}$  ni de ${y}$.

De l'hypoth\`ese de pincement, on  
 d\'eduit que l'entropie volumique $h_{\rm vol}(X)$ de $X$ est finie (avec $a(N-1)\leq h_{\rm vol}(X) \leq b(N-1)$ o\`u $N$ est la dimension de la vari\'et\'e \cite{GHL})  ;  il en est de m\^eme pour  l'exposant critique de $\Gamma$ puisque l'on a $\delta_\Gamma\leq h_{\rm vol}(X).$ 
 
 On peut aussi consid\'erer la fonction ``annulaire'' suivante :
 $$
  \forall {x, y}\in X, \forall R, \Delta >0 \quad n_\Gamma( R, \Delta):= \sharp\{\gamma \in \Gamma : R-\Delta \leq ({x}, \gamma\cdot{y})\leq R+\Delta\}.
 $$
(on \'ecrivant de m\^eme   $n_\Gamma(R, \Delta):= n_\Gamma({o}, {o}, R, \Delta)$). Comme ci-dessus, l'exposant $ \displaystyle \limsup_{R\to +\infty}{1\over R}\ln N_\Gamma({x}, {y}, R)$ ne d\'epend   ni de ${x}$  ni de ${y}$.

 Gr\^ace au lemme \'el\'ementaire suivant, nous avons en fait les \'egalit\'es 
 $$
 \limsup_{R\to +\infty}{1\over R}\ln N_\Gamma(  R)= \limsup_{R\to +\infty}{1\over R}\ln n_\Gamma( R, \Delta).
 $$

 \begin{lemma}\label{expocritique}
 Soit $(u_n)_{n \geq 0}$ une suite de r\'eels positifs telle que $\displaystyle \sum_{n\geq 0} u_n =+\infty$ 
   ; on pose $U_n:= u_0+\cdots + u_n$ pour tout $n\geq 0$. Alors, pour tout $s>0$, les s\'eries $\displaystyle \sum_{n\geq 0}u_n e^{-sn}$ et $\displaystyle \sum_{n\geq 0} U_n e^{-sn}$ convergent ou divergent simultan\'ement. En particulier, elles ont le m\^eme exposant critique $u$, donn\'e par
   $$u:= \limsup_{n\to +\infty} {1\over n}\ln u_n = \limsup_{n\to +\infty} {1\over n}\ln U_n \geq 0.$$
 \end{lemma} 
  L'exposant $\displaystyle \delta_\Gamma:=\limsup_{R\to +\infty}{1\over R}\ln N_\Gamma(  R)= \limsup_{R\to +\infty}{1\over R}\ln n_\Gamma( R, \Delta)$ est aussi appel\'e {\it exposant de Poincar\'e} de $\Gamma$, ou {\it exposant critique de $\Gamma$} ; en effet, la s\'erie dite {\it  de Poincar\'e } d\'efinie par
$$
  \forall {x, y}\in X, \forall s \geq 0  \quad {\mathcal P}_\Gamma( {x}, {y}, s):= \sum_{\gamma \in \Gamma}e^{-s({x}, \gamma\cdot{y})}
 $$ 
 diverge pour $s<\delta_{\Gamma}$ et converge pour  $s>\delta_{\Gamma}$.

 {\bf Exemples}
\begin{itemize}
\item Si $\Gamma=\langle h\rangle$ o\`u $h$ est une isom\'etrie hyperbolique, son exposant critique est nul (noter que la ``limsup'' est en fait une limite et que  la s\'erie de Poincar\'e diverge en $0$).
\item Dans $\mathbb H^2$, notons $p$ l'isom\'etrie parabolique $z\mapsto z+1$ ;  un calcul explicite en g\'eom\'etrie hyperbolique donne $(i, p^n\cdot i)= (i, i+n)=2\ln n+O(1)$ si bien que l'exposant     de  $\langle p\rangle$  vaut ${1\over 2}$. Noter que la ``limsup''  est aussi une limite ici et que  la s\'erie de Poincar\'e de $\langle p\rangle$ diverge en $1/2$. Remarquer aussi que si la courbure est constante \'egale \`a $-b^2$, l'exposant critique du groupe $\langle p\rangle$ devient alors \'egal \`a $b/2$. 
\item On peut modifier la m\'etrique hyperbolique de $\mathbb H^2$ de fa\c{c}on \`a avoir des larges bandes 2 \`a 2  disjointes $  \{z : a_n\leq  {\rm Im}z \leq b_{n}\}$, avec $b_n<a_{n+1}, $ o\`u la m\'etrique est alternativement \`a courbure $-a^2$ ou $-b^2$ (avec $b>a)$ Si la largeur des bandes $b_n-a_n$ cro\^it suffisamment vite, on aura alors 
$$
\liminf_{R\to +\infty}{1\over R}\ln N_{\langle p\rangle}(R)=a/2\quad \mbox{\rm} \quad < \quad  \limsup_{R\to +\infty}{1\over R}\ln N_{\langle p\rangle}(R)=b/2.$$
\item Si $\Gamma$ est co-compact dans $\mathbb H^d, d\geq 2$, son exposant critique vaut $d-1$ ; il en est de m\^eme en courbure variable, la valeur $d-1$ \'etant alors remplac\'ee par l'entropie volumique de $X$, d\'efinie par
\begin{equation}\label{entropievol}
h_{\rm vol}(X):= \limsup_{R\to +\infty}{1\over R}\ln {\rm Vol} (B_X({o}, R))
\end{equation}
et la ``limsup'' est une limite dans ce cas.  Cet argument provient  de \cite{Ma} o\`u il est d\'emontr\'e  que l'exposant critique d'un tel groupe co\"{i}ncide avec l'entropie topologique  du flot g\'eod\'esique sur le fibr\'e unitaire tangent de la vari\'et\'e compacte consid\'er\'ee (cet argument a \'et\'e repris par ailleurs par  G. Robert \cite{robert} pour \'etudier la croissance des groupes hyperboliques). Donnons quelques pr\'ecisions  : on choisit un domaine  fondamental $\mathcal D$ pour l'action de $\Gamma$, domaine que l'on peut supposer relativement compact  et dont on note $\Delta$ le diam\`etre, on a alors
$$
\bigcup_{\gamma\in \Gamma / (o, \gamma\cdot o)\leq R-\Delta}
\gamma\cdot \mathcal D
\subset B(o, R)\leq
\bigcup_{\gamma\in \Gamma / (o, \gamma\cdot o)
\subset
 R+\Delta}
\gamma\cdot \mathcal D.
$$ 
L'\'egalit\'e (\ref{entropievol}) s'en d\'eduit imm\' ediatement ; le fait que la ``limsup'' soit une limite  vient du fait que, si l'on pose $u_n:= n_\Gamma(n, 2\Delta)$,  on a, \`a une constante multiplicative pr\`es,
 $$
 u_{n+m}\leq u_n u_m.
 $$
 Il vient 
 \begin{equation}\label{conv-sousadd}
  (u_n)^{1/n} \to u:=\inf_{n\geq1}(u_n)^{1/n}= e^{\delta_\Gamma}
  \end{equation}
\item
 Si $\Gamma$ est un r\'eseau non uniforme  de $\mathbb H^d, d\geq 2$ (ie ${\rm vol}(\mathbb H^d/\Gamma)<+\infty$ mais $\mathbb H^d/\Gamma$ non compact), son exposant critique vaut $d-1$; il en est de m\^eme en courbure variable $1/4$-pinc\'ee (ie $ b^2\leq  4a^2)$ la valeur $d-1$ \'etant   remplac\'ee l\`a aussi par l'entropie volumique $h_{\rm vol}(X)$ de $X$ (voir \cite{DPPS1}$^($\footnote{l'id\'ee de base est que pour toute horoboule $\mathcal H$ et tout point ${x}$ sur l'horosph\`ere $\partial \mathcal H$,  l'ensemble $B({x}, R)\cap \mathcal H$ est approximativement \'egal \`a la  boule de rayon $R/2$ int\'erieure \`a $\mathcal H$ et tangente en ${x}$ \`a  $\partial \mathcal H$ ; son volume est donc comparable \`a celui d'une boule de rayon $R/2$ et n'influe donc pas sur la croissance du volume global lorsque la courbure est $1/4$-pinc\'ee}$^)$). Attention, lorsque la courbure n'est plus $1/4$-pinc\'ee, on peut avoir 
 $\delta_\Gamma < h_{\rm vol}(X)$.
\end{itemize}

Une premi\`ere question naturelle est de savoir si la ``limsup'' est une limite ou non ; la r\'eponse est ``non'' de fa\c{c}on g\'en\'erale, comme l'exemple du groupe parabolique ci-dessus le sugg\`ere. Cependant, d\`es que $\Gamma$ n'est pas \'el\'ementaire, nous avons le 
\begin{theorem} \label{exposantcritique}
(T. Roblin \cite{Rob1} \& D. Sullivan \cite{Su})Si $\Gamma$ est un groupe  non \'el\'ementaire, alors, pour tous ${x}$ et ${y} $ de $X$, on a 
$$\delta_\Gamma = \lim_{R\to +\infty}{1\over R}\ln N_\Gamma({x}, {y}, R).$$
De plus, il existe une constante $C:= C(\Gamma, {x}, {y}) >0$ telle que 
$$
N_\Gamma({x}, {y}, R)\leq C e^{\delta_\Gamma R}.
$$
\end{theorem}
La premi\`ere assertion de ce th\'eor\`eme  est due \`a Th. Roblin, tandis que la seconde d\'ecoule du c\'el\`ebre ``lemme de l'ombre '' de Sullivan ; la d\'emonstration de ces deux assertions repose sur l'existence de densit\'es conformes et la construction de Patterson de telles familles de mesures.
Nous proposons ici une approche radicalement diff\'erente et \'el\'ementaire (\cite{DPS}) qui permet en particulier de comprendre  certains ph\'enom\`enes qui apparaissent lorsque   $X$ est remplac\'e par un rev\^etement normal non simplement connexe.

\noindent D\'emonstration. On utilise de fa\c{c}on cruciale le lemme classique suivant  (voir par ex \cite{CI}) qui d\'ecoule du fait que la courbure est pinc\'ee :

\begin{lemma}\label{quasiiso}  Pour tout $\theta >0$,  il existe une constante  $D>0$ d\'ependant de $\theta$ et de la borne sup\'erieure de la courbure $-a^2$,   telle que pour tout triangle g\'eod\'esique
   $T$ de sommets ${x},{y},{z}\in X$ et dont l'angle en ${y}$ est plus grand que $\theta >0$, on a $$d({ x},{ z})\geq d({x},{ y})+d({y},{ z})-D.$$
 \end{lemma} 

  On fixe $\Delta >0$, $a, b >> \Delta$  et $\alpha, \beta \in \Gamma$ tels que 
  $a-\Delta \leq d({o}, \alpha \cdot {o})\leq a+\Delta$ et 
  $b-\Delta \leq d({o}, \beta \cdot {o})\leq b+\Delta$.
  
Le groupe $\Gamma$ \'etant non \'el\'ementaire, il poss\`ede un sous-groupe de type $\theta $-Schottky $H =\langle h_1, h_2\rangle $ avec $\theta>0$\ 
$^($\footnote{ le terme $\theta $-Schottky signifie  que, pour tous mots ``admissible'' $h$ et $h'$ en les lettres $h_1^{\pm 1}$ et $h_2^{\pm 1}$ et dont les premi\`eres lettres diff\`erent, l'angle en ${o}$ du triangle $ {o}, h\cdot{o}, h'\cdot{o}$ est $\geq \theta$}$^)$ ; on peut alors, selon la position relative des points $\alpha^{-1}\cdot {o}$ et $\beta\cdot{o}$,  associer de fa\c{c}on unique au couple $(\alpha, \beta)$ un couple $(\alpha', \beta)$, avec $\alpha' = \alpha h_i^{\pm 2}$, et  de telle sorte que  l'angle en {o} du triangle $(\alpha')^{-1}\cdot {o}, \ {o}, \  \beta\cdot{o}$ soit $\geq \theta_0>0$.
D'apr\`es le lemme pr\'ec\'edent,  on a approximativement
$$ ({o}, \alpha'\beta\cdot {o}) \simeq a+b. $$
Remarquons par ailleurs qu'il existe un entier $M\geq 1$ tel que tout \'el\'ement $\gamma \in \Gamma$ v\'erifiant $({o}, \gamma\cdot{o})\simeq a+b $  se d\'ecompose en au plus   $M$ fa\c cons distinctes de la forme $\gamma = \alpha'\beta$ avec  $({o}, \alpha'\cdot  {o})\simeq a  $ et $({o}, \beta \cdot {o})\simeq b.$ 

On prouve ainsi l'existence d'une constante $C>0$ et d'un entier $\kappa \geq 1$ tels que 
$$
\forall k, l >> 0 \quad u_ku_l\leq C \sum_{i=k+l-\kappa}^{k+l+\kappa} u_i 
$$
o\`u l'on a pos\'e $u_k:= n_\Gamma(k, \Delta)$  pour tout entier $k\geq 0$.
Lorsque la suite $(u_{k+1}/u_k)_{k \geq 1}$ est major\'ee$^($\footnote{cette propri\'et\'e est v\'erifi\'ee  en particulier lorsque $\Gamma$ est co-compact ou plus g\'en\'eralement convexe co-compact}$^)$, on voit, quitte \`a changer $u_l$ en $C'u_l$ avec $C'>0$,  que la suite $(u_k)_k$ est sur-multiplicative ; on sait alors que la suite $(u_k^{1/k})_k$ converge vers sa borne sup\'erieure. Le lemme s'en d\'eduit.

Plus g\'en\'eralement, on a le r\'esultat suivant qui permet de terminer  la d\'emonstration du th\'eor\`eme :  

\begin{fait} Soit $(u_n)_{n\geq 0}$ une suite de r\'eels positifs (non tous nuls) v\'erifiant
$$
\forall k, l \geq  \kappa \quad u_ku_l\leq  \sum_{i=k+l-\kappa}^{k+l+\kappa} u_i 
$$
o\`u $\kappa \geq 1$ est un entier   fix\'e 
et l'on pose $U_n:= u_0+\cdots+u_n$.

 La suite $(U_n^{1\over n})_{n\geq 0}$ converge alors vers un nombre $u\geq 1$ et il existe une constante $C>0$ telle que $\forall n\geq 0 \quad u_n\leq C u^n$.
\end{fait}
\fdem

 pose $U_k:= u_0+\cdots +u_k$ pour $k \geq 0$ et on obtient
$$
\forall k, l \geq 0 \quad u_k U_l\leq  C (2\kappa +1) U_{k+l+\kappa}.
$$
Quitte \`a multiplier chaque terme $u_k$ par une constante strictement positive, on peut \'ecrire
\begin{equation}\label{sur-mul1}
\forall k, l \geq 0 \quad u_k U_l\leq   U_{k+l+\kappa}.
\end{equation}
d'o\`u l'on d\'eduit, pour tout $q,b>0$ et $0\leq r<b$  
$$
U_{(b+\kappa)q+r+\kappa}\geq u_r U_{(b+\kappa)q}\geq u_r u_b U_{(b+\kappa)(q-1)}\geq ...
\geq  u_r (u_b)^q  
$$
et donc, en posant $m:=(b+\kappa)q+r+\kappa $
\begin{equation}\label{sur-mul2}
\Bigl(U_{m}\Bigr)^{1\over m}\geq (u_r)^{1\over m} (u_b)^{q\over m}. 
\end{equation}
On obtient, en faisant tendre $q$ puis $b$ vers $+\infty$
$$
\liminf_{+\infty} \Bigl(U_{m}\Bigr)^{1\over m}\geq 
\limsup_{b\to +\infty} \Bigl(u_{b}\Bigr)^{1\over b}:= u$$
avec $u\geq 1$ car $u_k\geq 1 i.s.$).
Comme $\displaystyle \limsup_{n\to +\infty} \Bigl(u_{n}\Bigr)^{1\over n}=\limsup_{n\to +\infty} \Bigl(U_{n}\Bigr)^{1\over n}$, on conclut que $\displaystyle u= \lim_{n\to +\infty} \Bigl(U_{n}\Bigr)^{1\over n}$ d'o\`u la premi\`ere assertion du Th\'eor\`eme.
L'in\'egalit\'e (\ref{sur-mul2}) nous donne alors, en faisant tendre $q$ vers $+\infty$
$$
\forall b \geq 0 \quad u\geq (u_b)^{1\over b+\kappa} 
$$
et la seconde assertion du Th\'eor\`eme s'en d\'eduit.\fdem
\noindent {\bf Question :\ }{\it  A-t-on aussi 
$\displaystyle \liminf_{R\to+\infty} (n_\Gamma(R, \Delta))^{1/R}= \limsup_{R\to+\infty} (n_\Gamma(R, \Delta))^{1/R}$ ? La r\'eponse semble \^etre ``Non'', il faudrait donc  caract\'eriser les groupes qui v\'erifient cette \'egalit\'e et donner des exemples explicites de groupes qui ne la v\'erifient pas.}

\section{Groupes divergents et trou critique}
 On dit qu'un groupe discret $\Gamma$ d'isom\'etries de  $X$  est {\bf divergent } lorsque sa s\'erie de Poincar\'e diverge en $s=\delta_\Gamma$ ; dans le cas contraire on dit que $\Gamma$ est {\bf convergent}.
 
 \noindent {\bf Exemples}
 \begin{itemize}
 \item L'exposant critique du groupe \'el\'ementaire  $\langle h\rangle$, o\`u $h$ est une isom\'etrie hyperbolique de $X$,  est nul ; ce groupe est donc divergent.
 \item L'exposant critique du groupe \'el\'ementaire $\langle p\rangle$, engendr\'e par l'isom\'etrie parabolique $z\mapsto z+1$ de $\mathbb H^2$   vaut $1/2$ puisque $(i, p^n\cdot i) = 2\ln n +O(1)$ ; de plus ce groupe est divergent.
 \item Il existe en courbure variable des groupes paraboliques convergents (voir \cite{DOP})
 \item Les groupes co-compacts ou convexe co-compacts sont divergents ; cela se d\'eduit d'un r\'esultat g\'en\'eral sur les mesures conformes, via encore une fois le lemme de l'ombre de Sullivan. On peut aussi le d\'emontrer  en
 utiliant (\ref{conv-sousadd}) qui entra\^ine imm\' ediatement $u_n\succeq e^{\delta_\Gamma n} 
 \quad \ ^($\footnote{o\`u la notation $a_n\succeq b_n$ signifie que la suite $(b_n\ / a_n)_n$ est  major\'ee}$^)$ et la divergence du groupe s'en d\'eduit imm\'ediatement.

  \end{itemize}
 
 \subsection{Un crit\`ere simple de s\'eparation des exposants}
 Le calcul explicite de l'exposant critique d'un groupe kleinien s'av\`ere \^etre une question souvent d\'elicate. Une question plus accessible est d\'ej\`a de comparer l'exposant critique d'un groupe avec celui de ses sous-groupes et de d\'egager un crit\`ere qui assure que l'on a une in\'egalit\'e stricte.

\begin{theorem}\cite{DOP} Soit $\Gamma$ un groupe kleinien et $H$ un sous-groupe de $\Gamma$ tel que 
\begin{enumerate}
\item $\Lambda_H\neq \Lambda_\Gamma$
\item $H$ divergent
\end{enumerate}
Alors, on a $\delta_H<\delta_\Gamma$.
\end{theorem}
D\'emonstration. On trouvera dans \cite{DOP} une d\'emonstration reposant sur l'existence de densit\'e $\delta_\Gamma$ conforme. Nous d\'eveloppons ici l'argument \'el\'ementaire suivant,  reposant sur la construction d'un sous-groupe de $\Gamma$ qui soit le produit libre de $H$ avec un sous-groupe \'el\'ementaire de $\Gamma$ engendr\'e par une isom\'etrie hyperbolique dont  les  points fixes n'appartiennent pas \`a $H$.

En effet, soit $\xi\in \Lambda_\Gamma\setminus \Lambda_H$ ; le point $\xi$ \'etant un point ordinaire de $H$, il existe un voisinage ouvert $U$ de $\xi$ tel que 
$U\cap \Lambda_H=\emptyset$.  L'action de $H$ sur son ensemble ordinaire \'etant propre et discontinue,   on peut supposer, quitte \`a r\'eduire $U$, que  
\begin{equation}\label{pingpong1}
\forall h\in H, h \neq id \qquad h(U)\subset \partial X\setminus U.
\end{equation}
Par ailleurs, l'action de $\Gamma$ \'etant minimale sur $\Lambda_\Gamma$, l'orbite d'un point  fixe d'une quelconque de ses isom\'etries est dense dans $\Lambda_\Gamma$, elle rencontre donc $U$. On en d\'eduit alors l'existence de $\gamma\in \Gamma\setminus H$ dont le point fixe attractif appartient \`a $U$ ; quitte \`a remplacer $\gamma$ par $\gamma^n$ avec $n$ grand, on peut supposer que $\gamma$ envoit l'ext\'erieur de $U$ dans $U$ (que $\gamma$ soit hyperbolique ou parabolique).

 Si on prend alors un  \'el\'ement hyperbolique $\gamma'\in \Gamma$ dont les points fixes sont distincts de ceux de $\gamma$,  et donc ext\'erieurs \`a $U$ quitte \`a r\'eduire $U$, l'isom\'etrie  $g:= \gamma^n\gamma'\gamma^{-n}$ a ses deux points fixes dans $U$, en rempla\c{c}ant si n\'ecessaire $g$ par une de ses puissances, on peut supposer que 
 \begin{equation}\label{pingpong2}
 \quad g^{\pm 1} (\partial X\setminus U)\subset U.
\end{equation}

Le sous-groupe $G$ de $\Gamma$ engendr\'e par $H$ et $g$ est alors un produit libre : il contient en particulier tous les mots de la forme
$h_1gh_2g\cdots h_kg$ avec $k\geq 1$ et $h_i\in H\setminus\{id\}$. On a alors, pour tout $s>\delta_H$
\begin{eqnarray*}
{\mathcal P}_\Gamma(s)&\geq &{\mathcal P}_\Gamma(s)\\
&\geq &\sum_{k\geq 1}\sum_{h_1, \cdots, h_k}e^{-s({o},h_1gh_2g\cdots h_k\cdot {o})} \\
&\geq &\sum_{k\geq 1}e^{-sk({o},g\cdot {o})}
\sum_{h_1, \cdots, h_k} e^{-s ({o},h_1\cdot {o})}\times \cdots \times e^{-s ({o},h_1\cdot {o})}\\
&\geq &\sum_{k\geq 1}\Bigl(e^{-s({o},g\cdot {o})}\times \sum_{h\in H^{*}} e^{-s ({o},h \cdot {o})}\Bigr)^k\\
&\geq &\sum_{k\geq 1}\Bigl(e^{-\delta_H({o},g\cdot {o})}\times \sum_{h\in H^{*}} e^{-s ({o},h \cdot {o})}\Bigr)^k.
\end{eqnarray*}
Le groupe $H$ \'etant divergent, on peut choisir $s_0>\delta_H$ suffisamment proche de $\delta_H$ de telle sorte que $e^{-\delta_H({o},g\cdot {o})}\times \sum_{h\in H^{*}} e^{-s_0 ({o},h \cdot {o})}>1$ si bien que ${\mathcal P}_\Gamma(s_0)=+\infty$ ; il vient $\delta_\Gamma \geq s_0>\delta_H$.

Attention ! les deux conditions  de ce th\'eor\`eme ne sont pas n\'ecessaires, mais l'on l'on en est pas tr\`es loin ! En effet
\begin{itemize}
\item Supposons $H$ convergent et $\Lambda_H\neq \Lambda _\Gamma$. On reprend la construction ci-dessus et l'on montre que le groupe $G_n= H*\langle g^n\rangle$ devient convergent lorsque $n$ est assez grand (utiliser le Lemme \ref{quasiiso}, voir \cite{DOP} pour la construction de groupes non \'el\'ementaires de et convergents).
\item Si le groupe $H$ est divergent et normal dans $\Gamma$, on a $\Lambda_H=\Lambda_\Gamma$ et $\delta_\Gamma= \delta_H$ (\cite{MY}, voir Th\'eor\`eme \ref{MY}).
\end{itemize}

La question de la convergence/divergence d'un groupe  Kleinien est int\'eressante \`a plus d'un titre; on renvoit le lecteur \`a \cite{Rob2} o\`u l'on trouve par exemple pour illustrer ce propos le th\'eor\`eme de dichotomie de Hopf, Tsuji \& Sullivan  avec une d\'emonstration   pr\'ecise.

Le th\'eor\`eme pr\'ec\'edent permet aussi d'exhiber de fa\c{c}on  simple des exemples de groupes   convergents. Consi\' erons par exemple un groupe de Schottky $\Gamma$ engendr\' e par deux isom\'etries hyperboliques $\alpha$ et $\beta$ et notons $G$ le sous-groupe de $\Gamma$ engendr\'e par la famille de  transformations $\{\alpha^{-n}\beta \alpha^n/ n\geq 0\}$. Les groupes $G$ et $H:=\alpha^{-1}G\alpha$ sont conjugu\'es, ils ont donc le m\^eme exposant critique $\delta$ et sont de m\^eme type (ie convergents ou divergents). Par ailleurs $\Lambda_G\subset \Lambda_H$ puisque $G\subset H$ mais  l'on a $\Lambda_G\neq \Lambda_H$ car  les points fixes de
$ \beta $ appartiennent \`a  $\Lambda_H$ mais pas \`a $\Lambda_G$. Les groupes $G$ et $H$ sont  donc convergents, d'apr\`es le th\'eor\`eme ci-dessus.

 \subsection{La question de la croissance forte}
 
  On consid\`ere un groupe Kleinien $\Gamma$ et $N$ un sous-groupe normal de $H$. Le groupe quotient $\bar \Gamma:= N\backslash \Gamma$ est le groupe des isom\'etries du rev\^etement Riemannien normal  $\bar X:= N\backslash X$  de $\Gamma \backslash X$, muni de la m\'etrique ``quotient''
  $$
  dist (\bar {x}, \bar {y}):=\inf \{(\bar {x}, n\cdot{y})/ n \in N\}.
  $$
  On note $\delta_{\bar \Gamma}$ l'exposant critique de $\bar \Gamma$ pour la m\'etrique correspondante ; on a $\delta_{\bar \Gamma}\leq \delta_{ \Gamma}$, et se pose alors de fa\c{c}on naturelle la question d'un crit\`ere assurant que cette in\'egalit\'e est stricte. Cette question peut \^etre vue comme la transposition dans un cadre Riemannien de celle de la croissance forte des groupes de type fini, pos\'ee par R. Grigorchuk et P. de la Harpe dans \cite{GH} 
  $^($\footnote{si $S$ est un ensemble fini de g\'en\'erateurs d'un groupe $G$, $\bar S$ l'ensemble correspondant dans $\bar G$ et si l'on munit $G$ de la m\'etrique des mots $\vert .\vert $  relativement \`a $S$, on pose 
  $\displaystyle w_{G,S}:= \lim_{n\to +\infty} \Bigl(\sharp\{g/\vert g\vert \leq n\}\Bigr)^{1\over n}$
  et l'on dit que 
 $G $ est \`a croissance forte relativement \`a $S$ si l'on a $w_{\bar G,\bar S}<w_{G,S}$. On montre ais\'ement que les groupes libres sont \`a croissance forte, il en est de m\^eme pour les produits amalgam\'es et les groupes hyperboliques \cite{AL} ; cette question est d'ailleurs \'etroitement  \`a celle de la croissance des sous-shift de type finie avec   mots interdits (voir par ex. \cite{HSW})
  }$^)$. 
 En reprenant les arguments d\'evelopp\'es ci-dessus, on peut montrer que   si $\bar \Gamma$ est divergent pour la m\'etrique quotient, alors on a $\delta_{\bar \Gamma}< \delta_\Gamma$ (\cite{DPPS2}).
 Il faut donc pouvoir pr\'eciser si $\bar \Gamma$ est divergent ou non.  \subsection{Sutr les groupes g\'eom\'etriquement finis}
 
 Nous avons vu dans le paragraphe pr\'ec\'edent des exemples de groupes divergents ; outre certains groupes \'el\'ementaires, les groupes co-compacts ou convexe co-compacts sont divergents.
 
 Curieusement, il est assez d\'elicat de d\'ecider si oui ou non un groupe non \'el\'ementaire $\Gamma$ est convergent et divergent. Une fois \'etudi\'ee la classe des  
 groupes co-compacts ou convexe co-compacts il est naturel de poser la m\^eme question pour les r\'eseaux non uniformes et leur g\'en\'eralisation naturelle : les groupes {\it g\'eom\'etriquement finis}.
 
 Introduisons quelques notations n\'ecessaires pour d\'efinir cette large classe de groupes kleiniens.On note $C (\Lambda_{\Gamma})$ l'enveloppe convexe de 
$\Lambda_{\Gamma}$; cet ensemble est invariant sous  l'action de 
$\Gamma$ et le quotient $N(\Gamma)= C (\Lambda_{\Gamma})/\Gamma$ est 
le coeur de Nielsen de la vari\'et\'e $  \Gamma\backslash X$. 
On notera $N_{\epsilon}(\Gamma)$ un $\epsilon$-voisinage
 de $N(\Gamma)$. Quand    l'ensemble   $N(\Gamma)$
est relativement compact ; on dit que $\Gamma$ est {\it convexe co-compact}. 
 On dit plus g\'en\'eralement que$\Gamma$ est g\'eom\'etriquement 
 fini lorsqu'il
 existe $\epsilon >0$ tel que
 $N_{\epsilon}$ soit de volume fini. 

La perte de compacit\'e de $N_{\epsilon}$ se traduit par 
l'appparition dans $\Lambda_{\Gamma}$ de points limites non coniques : les 
points   {\it paraboliques born\'ees}.

Rappelons qu'un point $\xi 
 \in \Lambda_{\Gamma}$ est dit  {\bf conique} (ou ``radial'') lorsque le segment g\'eod\'esique $[{o}, \xi)$ poss\`ede un voisinage qui contient une infinit\'e de points de l'orbite $\Gamma\cdot{o}$ et il est dit  {\bf parabolique born\'e} si son stabilisateur 
 $P$ est constitu\'e d'isom\'etries paraboliques et agit de fa\c{c}on relativement compacte sur   
 $\Lambda_{\Gamma}-\{\xi\}$.
 
 Quand $\Lambda_\Gamma$ ne contient que des points coniques, le groupe $\Gamma$ est convexe co-compact.
 
Rappelons que  la 
 finitude g\'eom\'etrique  peut \^etre caract\'eris\'ee de 
fa\c{c}on \'equivalente comme suit :
  
 -  pour tout $\epsilon >0$ le volume de $N_{\epsilon}(\Gamma)$ est fini.
 
 - pour tout $\epsilon >0$, la partie $\epsilon$-\'epaisse 
 $N(\Gamma)^{>\epsilon}$ est relativement compacte.
 
 - le groupe $\Gamma$ contient un nombre fini de classes de conjugaison  de 
 groupes paraboliques dont les points fixes  
 sont born\'es.
 
 - le coeur de Nielsen $N(\Gamma)$ peut se d\'ecomposer en $C_{0}\cup 
 C_{1} \cdots \cup C_{l}$ o\`u $C_0$ est un ensemble relativement 
 compact et o\`u, pour chaque $i = 1, \cdots, l,$
 il existe un groupe  parabolique $P_{i} \subset \Gamma$
et une 
 horoboule  ${\mathcal H}_{i}$ bas\'ee en $\xi_i$ tels que 
 $C_{i}$ soit isom\'etrique au 
 quotient de 
 ${\mathcal H}_i\cap C(\Lambda(\Gamma))$ par le groupe $P_{i}$
 (notons  que le point fixe $\xi_{i}$  de $P_{i}$ est alors 
 n\'ecessairement born\'e, que le groupe $P_{i}$ agit sur $C(\Lambda{G})\cap
 \partial \cH_{i}$ o\`u $\partial \cH_{i}$ d\'esigne l'horosph\`ere 
 qui
 borde l'horoboule  $ \cH_{i}$ et que cette action admet  un domaine 
 fondamental relativement compact).

Nous avons alors le
\begin{theorem} ( \cite{Su}, \cite{CI}\& \cite{DOP})
Soit $\Gamma$ un groupe g\'eom\'etriquement fini tel que, pour tout sous-groupe parabolique $\mathcal P$ de $\Gamma$ on ait  $\delta_\Gamma > \delta_{\mathcal P}$. 
Alors $\Gamma$ est divergent. 
\end{theorem}
L'approche initi\'ee par D. Sullivan et d\'evelopp\'ee jusqu'alors pour d\'emontrer ce r\'esultat repose sur l'existence d'une densit\'e conforme de dimension    $\delta_\Gamma$ et invariante par $\Gamma$ ; or, la d\'efinition m\^eme de la finitude g\'eom\'etrique, avec une partie \'epaisse relativement compacte, et une partie fine sur laquelle la dynamique de $\Gamma$ est tr\`es simple \`a contr\^oler, fait esp\'erer un autre argument de type sous-additivit\'e, comme dans le cas compact ou convexe co-compact.
 
 Pour ce faire, nous fixons  $\delta\in ]\delta^* _\mathcal P ; \delta_\Gamma[$ (o\`u  $\delta^* _\mathcal P$ d\'esigne le maximum des exposants critiques des sous-groupes paraboliques de $\Gamma$) et  nous introduisons la  quantit\'e  $v_\Gamma(R, \Delta)$ d\'efinie par 
$$v_\Gamma(R, \Delta) := e^{-\delta R}n_{   \Gamma}(    R, \Delta).$$
De part le choix de $\delta$, on a 
\begin{equation}\label{strict+}
\limsup_{R\to +\infty} {\ln v_\Gamma(R, \Delta)\over R} =\delta_\Gamma-\delta>0.
\end{equation}
Pour  tous r\'eels $a, b >0$ et $\Delta>1$ assez grand (sup\'erieur au diam\`etre de la partie \'epaisse $\mathcal C_0$), on a 
 \begin{equation}\label{sousadditif2} 
v_\Gamma(a+b, \Delta)\leq c \cdot  
\Bigl(\sum_{0\leq k\leq a} v_\Gamma(k, \Delta)\Bigr)  \cdot
\Bigl(\sum_{0\leq l\leq b} v_\Gamma(l, \Delta)\Bigr) 
\end{equation} 
o\`u $c$ est une constante $>0.^($\footnote{Pour d\'emontrer cette in\'egalit\'e, on fixe $\gamma$ tel que $\gamma\cdot{o}$ se trouve dans l'anneau $\{{x} / a+b-2\Delta\leq ({o}, {x})\leq  a+b+2\Delta\}$, on pose  $({o}, \gamma\cdot{o})=a+b+2\Lambda$  avec $-\Delta < \Lambda < \Delta$  et l'on note ${x}$ le point du segment g\'eod\'esique $[{o}, \gamma\cdot{o}]$ qui se trouve \`a distance $a+\Lambda $  de ${o}$  ; on somme alors sur $u\leq a$ et $v\leq b$, o\`u $u$  ( resp. $v$) repr\'esente la distance \`a parcourir \`a partir de ${x}$ sur le segment $[{x}, {o}]$ (resp. $[{x}, \gamma\cdot {o}]$  pour sortir de la zone mince de la vari\'et\'e quotient}$^)$.
 
 Posons  alors pour simplifier $w_n:= {v_\Gamma(n, \Delta)\over c}$, puis $W_n:= w_1+\cdots +w_n$ et $\tilde{W}_n:= 1+W_1+\cdots+W_n$. L'in\'egalit\'e
 $$
 \forall n, m \geq 1 \quad w_{n+m} \leq  W_n\times W_m
 $$
entra\^{i}ne successivement $W_{n+m}\leq W_n\tilde{W}_m$  puis $\tilde{V}_{n+m}\leq \tilde{V}_n\tilde{W}_m$.
Ainsi, la suite $(\ln \tilde W_n)_n$ est sous-additive et l'on a 
$$
\lim_{n\to +\infty}{\ln \tilde W_n\over n} = L:= \inf_{n \geq 1} {\ln \tilde W_n\over n}
$$
 En particulier $\tilde{W}_n\succeq e^{Ln}$ et la s\'erie $\sum_{n\geq 1} \tilde W_n e^{-Ln}$ diverge.
 
On  d\'eduit  alors du Lemme \ref {expocritique} que les  s\'eries $\sum_{n\geq 1} \tilde W_n e^{-sn}, \sum_{n\geq 1} W_n e^{-sn}$ et $\sum_{n\geq 1} w_n e^{-sn}$ admettent $L$ comme exposant critique. L'in\'egalit\'e  (\ref{strict+}) donne $L>0$ si bien que ces trois s\'eries divergent en $L$.  L'exposant critique de la s\'erie $\sum_{n\geq 1}  n_{   \Gamma}(    n, \Delta) e^{-sn}$ est  donc $\delta_\Gamma = L+\delta$ et cette s\'erie diverge en $\delta_\Gamma$.\fdem
 
\begin{remark} De ce qui pr\'ec\`ede, on d\'eduit que $\tilde{W_n}\succeq e^{Ln}$ et donc $W_n \succeq {e^{Ln}\over n}$. On peut en fait affiner cette derni\`ere in\'egalit\'e   en combinant l'estimation $\tilde{W_n}\succeq e^{Ln}$ avec l'in\'egalit\'e     
$W_n \preceq  e^{Ln}$ (cf Th\'eor\`eme  \ref{exposantcritique}
 ) : on peut ainsi obtenir de fa\c{c}on tr\`es \'el\'ementaire l'encadrement suivant 
 $$
 e^{\delta_\Gamma R}\preceq N_\Gamma(R)\preceq e^{\delta_\Gamma R}.
 $$
\end{remark}

\  \section{Mesure de Patterson et exposant critique}
  On appelle {\it densit\'e}  sur $\partial X$ une famille  $\mu=(\mu_x)_{x \in X}$ de mesures positives finies sur $\partial X$. Une telle densit\'e est dite {\it conforme de dimension $\delta\geq 0$} lorsque pour tous $x, x' \in X$, la mesure $\mu_{x'}$ est absolument continue  par rapport \`a $\mu_x$ et que sa d\'eriv\'ee de Radon-Nikodym est donn\'ee par la formule
  $$
  {d\mu_{x'}\over d\mu_x}(\xi) = e^{-\delta({\mathcal B}_\xi(x', x)}
  $$
  o\`u ${\mathcal B}_\xi(x', x)$ d\'esigne la  fonction de Busemann en $\xi$ d\'efinie par 
  $$
  {\mathcal B}_\xi(x', x):= \lim_{z \to \xi} (x', z)-(x, z).
  $$
  Cette densit\'e est dite {\it invariante par $\Gamma$} si pour tout $\gamma \in \Gamma$ et tout $x \in X$ on a 
  $$
  \gamma^*\mu_{\gamma\cdot x} = \mu_x   $$
  o\`u la mesure $\gamma^*\mu_x$ est d\'efinie par $\gamma^*\mu_x(A)= \mu_x(\gamma A)$ pour tout bor\'elien $A$ de $\partial X$.

  On  note  $Conf(\Gamma, \delta)$ l'ensemble des  densit\'es conformes $\Gamma$-invariantes et de dimension $\delta$, avec la condition de normalisation $\Vert
  \mu_o\ \Vert = 1$.
  
  \noindent Rappelons le proc\'ed\'e de Patterson permettant de 
montrer que $Conf(\Gamma, \delta_\Gamma)$ est non vide.
Pour chaque $s>\delta_{\Gamma}$ et chaque point $x \in X$
 on note $\mu_{x, y}^s$ la mesure orbitale
 $$
\mu_{x, y}^s= \frac{1}{{\mathcal P}_{\Gamma}(o, y, s)}
\sum_{\gamma \in \Gamma}\exp (-s  (x, \gamma\cdot y)) \epsilon_{\gamma\cdot y}$$
o\`u $\epsilon_{\gamma\cdot o}$ d\'esigne la masse de Dirac en $\gamma.y$.
Lorsque le groupe $\Gamma$ est de  type divergent,
toute valeur d'adh\'erence (pour la topologie de la convergence 
\'etroite)
de la famille 
$(\mu_{x, y}^s)_{x, y, s}$ est port\'ee par
$\Lambda_{\Gamma}$; on peut alors  montrer que lorsque $s \to \delta_{\Gamma}$ 
par valeurs sup\'erieures
 la famille de mesures $(\mu_{x, y}^s)_{s}$  converge \'etroitement 
 vers une  mesure $\mu_{x, y}$
port\'ee par $\Lambda_{\Gamma}$ et v\'erifiant les deux conditions suivantes
$$
\mu_{x',y}(.) = \exp (-\delta{\mathcal B}_{.}(x', x))\mu_{x, y}(.)
\quad \mbox{\rm et} \quad g^{*}\mu_{x, y} =\mu_{g^{-1}\cdot x, y} 
$$
o\`u $g^{*}\mu_{x, y} $ est la mesure sur $\partial X$ d\'efinie par 
$g^{*}\mu_{x, y}(B) = \mu_{x, y}(g B)$ pour tout bor\'elien $B $ 
de $\partial X$. 
On dit que la famille $(\mu_{x, y})_{x \in X}$ est une 
{\it densit\'e $\Gamma$-conforme d'exposant $\delta_{\Gamma}$}.

\noindent On a vu qu'il peut \^etre d\'elicat  de montrer qu'un groupe $\Gamma$ est 
de type divergent ;   D. Sullivan a \'etabli cette 
propri\'et\'e pour les groupes g\'eom\'etriquement finis, en 
\'etudiant le type des densit\'es $\delta_{\Gamma}$-conformes de ces groupes. 
Pour ce faire, il faut pouvoir construire de telles densit\'es  
$\delta_{\Gamma}$-conformes,  sans savoir \`a priori si $\Gamma$ est de type 
convergent ou 
divergent  ;  en
utilisant un argument du \`a Patterson, on modifie
l\'eg\`erement la s\'erie de Poincar\'e en posant
$${\mathcal P}'_{\Gamma}(x, y, s) = \sum_{g \in \Gamma}e^{-s(x, \gamma\cdot y)}h(d(x, \gamma\cdot y))$$
o\`u $h$ est une fonction croissante de $\R^{+}$ dans $\R^{+}$ telle que
les s\'eries ${\mathcal P}_{\Gamma}(x, y, s) $ et ${\mathcal P}'_{\Gamma}(x, y, s)$ aient le m\^eme 
exposant critique et
$$
\forall \eta >0,  \exists  t_{\eta}>0,  \forall t \geq t_{\eta},  \forall s 
\geq 0 \quad
 h(t+s) \leq h(t) e^{\eta s}.
$$

\noindent  Int\'eressons nous maintenant aux propri\'et\'es locales 
 d'une densit\'e $\Gamma$-conforme. 
Nous avons le 

\begin{lemma} (-Th\'eor\`eme de l'ombre de Sullivan-)  
Soit $\Gamma$ un groupe non \'el\'ementaire et $\mu $ une densit\'e 
$\Gamma$-conforme d'exposant $\alpha$. Il existe $C>1$ et $r_{0}>0$ tel 
que pour tout $r\geq r_{0}$ et tout $g \in \Gamma$ on ait 
$$
\frac{1}{C} e^{-\alpha d(o, \gamma\cdot o)} \leq  
\mu_{x, y}(O_{x}(\gamma\cdot o, r))\leq C  e^{-\alpha d(o, \gamma\cdot o)+2\alpha }.
$$
\end{lemma}

Soulignons  en cependant 
quelques  cons\'equences  importantes    :
\begin{itemize}
\item  en remarquant que pour tout $\Delta >0$,  la famille des ombres 
$O_{x}(\gamma\cdot o, r)$  avec $\gamma \in \Gamma$ et $R-\Delta \leq (o, \gamma\cdot o)  \leq R+\Delta$ forme un recouvrement  \`a multiplicit\'e uniform\'ement born\'ee par rapport  \`a $R$, on montre que  l'existence d'une
densit\'e $\Gamma$-conforme $\mu$ d'exposant $\alpha$ 
entra\^ine 
$$
\sharp\{ \gamma \in \Gamma/  (o, \gamma\cdot o) \leq R\} \preceq e^{\alpha R}
$$
et donc  $\alpha \geq \delta_{\Gamma}$ par d\'efinition de 
l'exposant  de Poincar\'e de $\Gamma$. Il n'existe donc pas de densit\'e 
$\Gamma$-conforme d'exposant  $< \delta_{\Gamma}$.
\item  l'ensemble radial de $\Gamma$ est \'egal  \`a
$\displaystyle \Lambda_\Gamma^{rad}:= \cup_{r >0}\limsup_{\stackrel{\gamma \in \Gamma}{\gamma \to \infty}} O_{ o}(\gamma\cdot o, r)
 $ ; par cons\'equent, 
si $\displaystyle {\sum_{n\geq 1}e^{-\alpha d(o, \gamma\cdot o)}<+\infty }$
(ce qui est le cas lorsque $\alpha >\delta_\Gamma$) alors toute densit\'e
$\alpha$-conforme donne une mesure nulle \`a $\Lambda_\Gamma^{rad}$.
\item  lorsque $\Gamma$ est divergent, on montre que $\mu_{x,y}( \Lambda_\Gamma^{rad})>0$ et que cette propri\'et\'e entra\^{i}ne l'ergodicit\'e de l'action de $\Gamma$ relativement \`a la famille $(\mu_{x, y})_{x, y}$. Ainsi, dans ce cas, il existe,  \`a une constante multiplicative pr\`es,  une unique densit\'e $\delta_\Gamma$-conforme ; en particulier,  chaque  mesure $\mu_{x, y}$ construite ci-dessus ne d\'epend   ni du point $y$ ni de la sous-suite $(s_k)_k$ qui a permis de la construire.  
\item Lorsque $\Gamma$ est co-compact ou plus g\'en\'eralement  convexe co-compact, on peut montrer que  $Conf(\Gamma, \delta)\neq \emptyset$ si et seulement si $\delta_\Gamma$ ;  dans le cas contraire, $Conf(\Gamma, \delta)\neq \emptyset$ pour tout $\delta\geq \delta_\Gamma$. 

\end{itemize}

Remarquons que pour tout \'el\'ement $g$ du normalisateur $ N(\Gamma)$ de $\Gamma$ dans le groupe des isom\'etries de $X$ et toute densit\'e  $\mu\in Conf(\Gamma, \delta)$, la famille de mesures
$(\nu_x^g)_{x\in X}$  d\'efinie par 
\begin{equation}\label{nu_gconforme}\nu_x^g:= {1\over\Vert \mu_{g\cdot o}\Vert}g^*\mu_{g\cdot x}
\end{equation} appartient aussi  \`a $Conf(\Gamma, \delta)^($\footnote{il suffit de v\'erifier que $\nu^g$ est invariante par $\Gamma$ ;  pour tout $\gamma \in \Gamma$,  on a $g\gamma = \tilde{\gamma} g$  avec $\tilde \gamma = g\gamma g^{-1} \in \Gamma$ 
 ; le fait que $\mu$ est $\Gamma$-invariante entra\^{i}ne la suite d'\'egalit\'es suivantes
$$
\gamma^*\nu_{\gamma\cdot x}^g=\gamma^*g^*\mu_{g\gamma\cdot x}= g^*\tilde{\gamma}^*\mu_{g\gamma\cdot x}= 
g^*\mu_{\tilde{\gamma}^{-1}g\gamma\cdot x}=
g^*\mu_{ g\cdot  x}=
\nu_{x}^g.
$$
}$^)$

Supposons  $\Gamma$  divergent. L'ensemble $Conf(\Gamma, \delta_\Gamma)$ est alors r\'eduit  \`a un point (not\'e   $(\mu_x)_{x \in X}$) et  l'on a  $(\nu_x^g)_x=(\mu_x)_x$ pour  tout $g \in N(\Gamma)$. 
%il existe donc  une constante $c(g)>0$ telle que $\nu_x^g= c(g)\mu_x.$ 
%Notons que 
%$$\Vert \mu_x\Vert= \lim_{s\to\delta_\Gamma^+}{{\mathcal P}_\Gamma(x, o,s)\over {\mathcal P}_\Gamma(o, o,s)}.
%$$
 %L'\'egalit\'e $\nu_x^g= c(g)\mu_x$ s'\'ecrit  $g^*\mu_{g\cdot x}=c(g)\Vert \mu_{g\cdot o}\Vert \mu_x$ et donc 
  %$\displaystyle c(g) \Vert \mu_{g\cdot o}\Vert= {\Vert \mu_{g\cdot x}\Vert   \over \Vert \mu_x\Vert }$ pour tout $g \in N(\Gamma)$ et $x \in X$.
  %Il vient d'une part
  %$$
  %c(g)= \Vert \mu_{g\cdot o}\Vert = \lim_{s\to\delta_\Gamma^+} {{\mathcal P}_\Gamma(g\cdot o, o,s)\over {\mathcal P}_\Gamma(o, o,s)}
  %$$
  %et d'autre part
  %$$
  %{1\over c(g)}= \Vert \mu_{g^{-1}\cdot o}\Vert =\lim_{s\to\delta_\Gamma^+} {{\mathcal P}_\Gamma(g^{-1}\cdot o, o,s)\over {\mathcal P}_\Gamma(o, o,s)}.
  %$$
%L'\'egalit\'e ${\mathcal P}_\Gamma(x, y, s)={\mathcal P}_\Gamma( y, x, s)= {\mathcal P}_\Gamma(g\cdot x, g\cdot y, s)$ pour tout $x, y \in X$ et $g \in N(\Gamma)$ permet de conclure $c(g)= {1\over c(g)}=1$. 
Ainsi $Conf(N(\Gamma), \delta_\Gamma)\neq \emptyset$. Par cons\'equent   $\delta_\Gamma \geq \delta_{N(\Gamma)}$ et donc en fait $\delta_\Gamma = \delta_{N(\Gamma)}$.

On a ainsi montr\'e le 
\begin{theorem} \label{MY} \cite{MY}
Soit $\Gamma$ un groupe kleinien divergent et $G$ un groupe kleinien qui contient $\Gamma$ comme sous-groupe normal. Alors $\delta_\Gamma=\delta_G$ (et $G$ est aussi divergent).
\end{theorem}

Plus g\'en\'eralement, Th. Roblin  \`a introduit la notion de {\bf principe des ombres}  d\'efinie comme suit :
{\it 
On dit que le sous-ensemble $Y \subset C(\Lambda_\Gamma)\subset X$ invariant sous l'action de $\Gamma$ v\'erifie le principe des ombres en $\delta \geq \delta_\Gamma$ si pour toute densit\'e $\mu\in  Conf(\Gamma, \delta)$ et tous $x, y \in Y$ on a
$$
{1\over C}\Vert\mu_y\Vert e^{-\delta(x, y)}\leq \mu_x(O_x(y, r))\leq C
\Vert\mu_y\Vert e^{-\delta(x, y)}
$$
o\`u $r$ et $C$ sont des constantes positives ne d\'ependant que  de $Y$ et de $\delta$ (et non de $x, y$ et $\mu$).
}

D'apr\`es le Lemme de l'ombre de Sulivan, l'orbite $\Gamma\cdot o$ satisfait au principe des ombres ; plus g\'en\'eralement on peut aussi consid\'erer l'orbite de $o$ sous l'action du normalisateur de $\gamma$ puisque l'on a le
\begin{theorem}
Soit $N(\Gamma)$ le normalisateur de $\Gamma$ dans le groupe des isom\'etries de $X$. Alors, l'ensemble $N(\Gamma)\cdot o$ v\'erifie le principe des ombres pour tout $\delta\geq \delta_\Gamma$.
\end{theorem}
La  d\'emonstration est analogue  \`a celle du lemme de Sullivan ; on utilise de fa\c{c}on cruciale le fait que pour tout $\mu\in Conf(\Gamma, \delta)$ et tout $g\in N(\Gamma)$ la famille $(\nu^g)_{x\in X}$ d\'efinie par la formule (\ref{nu_gconforme}) est une densit\'e $\delta$-conforme.

Le fait qu'un ensemble $\Gamma$-invariant $Y$ satisfasse au principe des ombres entra\^{i}ne que la suite $(\Vert \mu_y\Vert)_{y\in Y}$ cro\^{i}t de fa\c{c}on contr\^ol\'ee. En effet, on a le 
\begin{theorem}
Soit $N(\Gamma)$ le normalisateur de $\Gamma$ dans le groupe des isom\'etries de $X$. Alors, pour toute densit\'e $\mu \in Conf(\Gamma, \delta)$, l'exposant critique de la s\'erie 
$\displaystyle 
\sum_{g\in N(\Gamma)}\Vert \mu_{g\cdot o}\Vert e^{- s (o, g\cdot o)} 
 $ est inf\'erieur ou \'egal \`a $\delta$.
\end{theorem}
qui admet de fa\c{c}on imm\'ediate le corollaire suivant
\begin{corollary}
Pour tout groupe kleinien $\Gamma$  de $X$ on a $\displaystyle \delta_{\Gamma}\geq {1\over 2}\delta_{N(\Gamma)}.$

De plus, si $N(\Gamma)$ est divergent alors cette in\'egalit\'e est stricte.
\end{corollary}

Posons pour simplifier $G:= N(\Gamma)$. La premi\`ere  assertion de ce corollaire est une cons\'equence directe du Th\'eor\`eme pr\'ec\'edent. Pour la seconde assertion, on  raisonne par l'absurde, on suppose que $\delta=\delta_\Gamma= \delta_{G}/2$, on note $\nu$ l'unique  \'el\'ement de $Conf( G,\delta_{G}) $ et l'on consid\`ere un \'el\'ement $\mu$ de $Conf( \Gamma,\delta_{ \Gamma}) $ que l on suppose ergodique.  Puisque $G$ est divergent, $\nu_o$ est port\'ee par $\Lambda_G^{rad}(r)$ pour $r$ assez grand . Si $B$ est un bor\'elien  de $\partial X$, on choisit un ouvert $V$ le contenant et on note $Z$ l'ensemble des points $z$ de $G\cdot o$ tel que $O_o(z, r)\subset V$ ; par un argument de type recouvrement de Vitali,  on peut extraire de $Z$ une sous famille finie $Z^*$ telle que  les ombres  correspondantes $O_o(z, r)$ soient deux \`a deux disjointes et
$\displaystyle \bigcup_{z\in Z}O_o(z, r)\subset \bigcup_{z\in Z^*}O_o(z,5 r). $
On a alors
\begin{eqnarray*}
\nu_o(B)=\nu_o(B\cap \Lambda^{rad}_G(r))&\leq& \sum_{z\in Z^*}\nu_o( O_o(z,5 r))
\\
&\preceq&\sum_{z\in Z^*}e^{-\delta_G (o,z)}\\
&=&\sum_{z\in Z^*}e^{-2 \delta_\Gamma (o,z)}\\
&\preceq&\sum_{z\in Z^*}\Vert \mu_y\Vert e^{-\delta_\Gamma (o,z)}\\
&\preceq&\sum_{z\in Z^*}\mu_o( O_o(z, r))\\
&\preceq& \mu_o(V)
\end{eqnarray*}
Il vient $\nu_o(B)\preceq \mu_o(B)$, l'ouvert $V$ \'etant arbitraire, et donc $\nu_o=\phi(.)\mu_o$ pour une certaine densit\'e $\phi\geq 0$. Or, pour tout $\gamma \in \Gamma$ on a $$\phi(\gamma \cdot \xi) = \phi(\xi)
 e^{-(\delta_G-\delta_\Gamma)\mathcal B(\gamma\cdot o, o)}\quad \mu_{ o}(d\xi)\quad p.s.
 ^(\footnote{
pour toute fonction bor\'elienne positive $\psi$ sur $\partial X$  on a  
d'une part
$$
\int\psi(\gamma^{-1} \xi)\nu_o(d\xi)
= \int\psi(  \xi)\nu_{\gamma\cdot o}(d\xi)=
\int\psi(  \xi)e^{-\delta_G\mathcal B(\gamma\cdot o, o)}\nu_{o}(d\xi)
=\int\psi(  \xi)e^{-\delta_G\mathcal B(\gamma\cdot o, o)}\phi(\xi)\mu_o(d\xi)
$$
et d'autre part
$$
\int\psi(\gamma^{-1} \xi)\nu_o(d\xi)= \int\psi(\gamma^{-1} \xi)\phi(\xi)\mu_o(d\xi)
= \int\psi(  \xi)\phi(\gamma\cdot \xi)e^{-\delta_\Gamma\mathcal B(\gamma\cdot o, o)}\mu_{ o}(d\xi)$$
d'o\`u
$\phi(\gamma\cdot \xi)
= \phi(\xi)
 e^{-(\delta_G-\delta_\Gamma)\mathcal B(\gamma\cdot o, o)}\quad \mu_{ o}(d\xi)$p.s.
}^)$$
La mesure $\mu_o$ \'etant ergodique, la fonction quasi-invariante $\phi$ est $\mu_o$-p.s. strictement positive ; les mesures $\mu_o$  et $\nu_o$ sont donc  \'equivalentes et par cons\'equent  quasi-invariantes par $G$. L'ergodicit\'e de $\mu_o$ montre que $\mu$ est   construite \`a partir d'un caract\`ere de $G/\Gamma$ ce qui entra\^{i}ne de facto $\delta\geq \delta_{G}$. Contradiction.\fdem

On a donc ${\delta_{N(H)}\over 2}\leq \delta_\Gamma \leq \delta_{N(H)}$, la premi\`ere in\' egalit\'e \'etant stricte d\`es que $G$ est divergent. Par ailleurs, on a vu que si c'est $\Gamma$ qui est divergent alors $\delta_\Gamma=\delta_G$.

Ainsi, si $\delta_\Gamma \in [{\delta_G\over 2}, \delta_G[$, le groupe $\Gamma$ est convergent   ; par contre lorsque $\delta_\Gamma=\delta_G$, on peut avoir $\Gamma$ convergent ou divergent (voir les rev\^etements ab\'eliens de vari\'et\'es compactes).

On a aussi le 
\begin{theorem} Si $N(\Gamma)/\Gamma$ est moyennable, alors $\delta_\Gamma =\delta_{N(H)}$.
\end{theorem}
La r\'eciproque de ce th\'eor\`eme est partiellement d\'emontr\'ee lorsque $G:= N(\Gamma)$ est convexe co-compact (travail de R. Brooks, dans $\mathbb H^{n+1}$ avec la condition $\delta_{G}>n/2$ \cite{Br})  ; la question de savoir si la r\'eciproque est vraie pour tous les groupes convexe co-compact, et plus g\'en\'eralement pour les groupes g\'eom\'etriquement finis, reste ouverte actuellement.

\end{document}